\documentclass[11pt]{article}
\usepackage{amsfonts}
\usepackage{amsmath}
\usepackage{amsthm}
\usepackage{url}
\usepackage[all]{xy}
\usepackage{graphicx}
\usepackage{color}
\usepackage{url}

\long\def\symbolfootnote[#1]#2{\begingroup\def\thefootnote{\fnsymbol{footnote}}
\footnote[#1]{#2}\endgroup}

\def\Mod{{\mbox{\rm {Mod}}}}
\def\Aut{{\mbox{\rm {Aut}}}}
\def\Out{{\mbox{\rm {Out}}}}
\def\Hom{{\mbox{\rm {Hom}}}}
\def\Inn{{\mbox{\rm {Inn}}}}
\def\Homeo{{\mbox{\rm {Homeo}}}}

\def\PSL{{\mbox{\rm {PSL}}}}
\def\PML{{\mathbb{P} \mathcal{ML}}}

\newcommand{\BF}{\mathbb F}

\newcommand{\CT}{\mathcal T}

\newcommand{\CZ}{\mathcal Z}
\newcommand{\CS}{\mathcal S}
\newcommand{\CM}{\mathcal M}

\newtheorem{conjecture}{Conjecture}

\newtheorem{corollary}[conjecture]{Corollary}

\newtheorem{lemma}[conjecture]{Lemma}

\newtheorem{proposition}[conjecture]{Proposition}
\newtheorem{question}{Question}
\newtheorem*{remark}{Remark}
\newtheorem*{definit}{Definition}
\newtheorem{theorem}[conjecture]{Theorem}

\newtheorem*{namedtheorem}{\theoremname}
\newcommand{\theoremname}{testing}
\newenvironment{named}[1]{\renewcommand{\theoremname}{#1}\begin{namedtheorem}}{\end{namedtheorem}}

\title{Injections of mapping class groups.}
\author{Javier Aramayona\thanks{Partially supported by M.E.C. grant MTM2006/14688}, Christopher J. Leininger\thanks{Partially supported by NSF grant DMS-0603881}  $\,$ and Juan Souto\thanks{Partially supported by NSF grant DMS-0706878 and Alfred P. Sloan Foundation}}

\begin{document}

\maketitle

\begin{abstract}
We construct new monomorphisms between mapping class groups of surfaces.  The first family of examples injects the mapping class group of a closed surface into that of a different closed surface.  The second family of examples are defined on mapping class groups of once-punctured surfaces and have quite curious behaviour. For instance, some pseudo-Anosov elements are mapped to multi-twists. Neither of these two types of phenomena were previously known to be possible although the constructions are elementary.
\end{abstract}

\section{Introduction}

Let $\Sigma_{g,k}$ be the closed orientable surface of genus $g$ with $k$ marked points and $\Homeo(\Sigma_{g,k})$ the group of homeomorphisms of $\Sigma_{g,k}$ which map the set of marked points to itself. If we endow $\Homeo(\Sigma_{g,k})$ with the compact open topology then $\Homeo_0(\Sigma_{g,k})$, the connected component of the identity, is the normal subgroup consisting of those elements which are isotopic to the identity relative to the marked points. The (extended) mapping class group of $\Sigma_{g,k}$ is the quotient group
$$\Mod(\Sigma_{g,k})=\Homeo(\Sigma_{g,k})/\Homeo_0(\Sigma_{g,k}).$$
Throughout this note we will assume without further mention that either $g\ge 2$, or $g=1$ and $k\ge 2$, or $g=0$ and $k\ge 5$.

The problem of studying homomorphisms between mapping class groups has received considerable attention. Automorphisms of mapping class groups were classified by Ivanov \cite{ivanovauto}, Korkmaz \cite{korkmazauto} and Luo \cite{luoauto}; and injections from finite index subgroups of a mapping class group to itself by Ivanov \cite{ivanovinjective}, Irmak \cite{irmakI,irmakII}, Behrstock-Margalit \cite{behrstockmargalit} and Shackleton \cite{shackleton}. Recently Harvey-Korkmaz \cite{harveykorkmaz} proved that there are no non-trivial homomorphisms $\Mod(\Sigma_{g,0})\to\Mod(\Sigma_{h,0})$ with $h<g$ and $g\ge 3$.

Besides the Harvey-Korkmaz theorem, the situation for homomorphisms between distinct mapping class groups is not well-understood and much of the research has focused on injective homomorphisms. For example, the fact that every finite subgroup of $\Mod(\Sigma_{g',1})$ is cyclic implies that there are no injective homomorphisms $\Mod(\Sigma_{g,0})\to\Mod(\Sigma_{g',1})$; see \cite{Mess}.  A stronger result is due to Ivanov-McCarthy \cite{ivanovmccarthy}, who proved that there are no injective homomorphisms $\Mod(\Sigma_{g,k})\to\Mod(\Sigma_{g',k'})$ for all but finitely many pairs $(g,k)$ and $(g', k')$ such that $3g'+k'-3g-k=1$.  However, in the same paper they observed that if $k>0$ then every characteristic cover $\Sigma_{g',k'}\to\Sigma_{g,k}$ induces an injective homomorphism $\Mod(\Sigma_{g,k})\to\Mod(\Sigma_{g',k'})$; see Section \ref{sec:hom cover}. Previously, Birman-Hilden \cite{Birman-Hilden} showed that a hyperelliptic branched cover $\Sigma_{g,0}\to\Sigma_{0,2g+2}$ also induces an injective homomorphism $\Mod(\Sigma_{0,2g+2})\to\Mod(\Sigma_{g,0})$. In both constructions, the marked points play a crucial role. The first goal of this note is to construct examples of injective homomorphisms between mapping class groups of surfaces without marked points.

\begin{theorem} \label{main}
For every $g\ge 2$ there is $g'>g$ and an injective homomorphism
\[\phi:\Mod(\Sigma_{g,0}) \to \Mod(\Sigma_{g',0}). \]
\end{theorem}

As in Ivanov-McCarthy \cite{ivanovmccarthy}, the idea of the proof of Theorem \ref{main} is to construct a finite-sheeted covering $\Sigma_{g',0}\to\Sigma_{g,0}$ with the property that every $f\in\Homeo(\Sigma_{g,0})$ has a distinguished lift. In order to do so we combine elementary covering theory and basic finite group theory such as the Sylow theorems.

Once the existence of injective homomorphism between mapping class groups has been sufficiently established, it is an interesting problem to study how these homomorphisms arise. Let $\phi:\Mod(\Sigma_{g,k})\to\Mod(\Sigma_{g',k'})$ be an injective homomorphism. Denoting by $\CT(\Sigma_{g,k})$ and $\CT(\Sigma_{g',k'})$ the Teichm\"uller spaces of $\Sigma_{g,k}$ and $\Sigma_{g',k'}$, it follows from the solution of the Nielsen realization problem \cite{kerckhoff} that there is a continuous $\phi$-equivariant map 
$$\Phi:\CT(\Sigma_{g,k})\to\CT(\Sigma_{g',k'}).$$
If the homomorphism $\phi$ is constructed as in \cite{Birman-Hilden}, \cite{ivanovmccarthy} or Theorem \ref{main} then it is not difficult to see that, endowing both Teichm\"uller spaces with the Teichm\"uller metric, $\Phi$ can be chosen to be an isometric embedding. In particular, any such $\phi$ is type-preserving in the sense that it maps pseudo-Anosov and reducible elements of $\Mod(\Sigma_{g,k})$ to again pseudo-Anosov and reducible elements of $\Mod(\Sigma_{g',k'})$. Moreover, the isometric embedding $\Phi$ extends continuously to an embedding of the Thurston compactification $\bar\CT(\Sigma_{g,k})=\CT(\Sigma_{g,k})\cup \PML(\Sigma_{g,k})$ of $\CT(\Sigma_{g,k})$ to the Thurston compactification $\bar \CT(\Sigma_{g',k'})$ of $\CT(\Sigma_{g',k'})$. Here $\PML(\Sigma_{g,k})$ is the space of projective measured laminations on $\Sigma_{g,k}$. In particular, the limit set $\Lambda_{\phi(\Mod(\Sigma_{g,k}))}$, in the sense of McCarthy-Papadopoulos \cite{McCarthy-Papadopoulos}, of the image of $\phi$ is a proper subset of $\PML(\Sigma_{g',k'})$ homeomorphic to $\PML(\Sigma_{g,k})$. 

Our next goal is to construct injective homomorphisms from $\Mod(\Sigma_{g,1})$ into some other mapping class group for which these properties fail.

\begin{theorem}\label{non-geometric}
For every $g$ there is $g'$ and an injective homomorphism $\alpha:\Mod(\Sigma_{g,1})\to\Mod(\Sigma_{g',1})$ with the following properties:
\begin{enumerate}
\item There are pseudo-Anosov mapping classes in $\Mod(\Sigma_{g,1})$ whose image under $\alpha$ is a multi-twist. In particular, $\alpha$ is not type-preserving.
\item The limit set $\Lambda_{\alpha(\Mod(\Sigma_{g,1}))}$ of the image of $\alpha$ is the whole space $\PML(\Sigma_{g',1})$ of projective measured laminations on $\Sigma_{g',1}$. 
\item In particular, there is no $\alpha(\Mod(\Sigma_{g,1}))$-invariant subset of $\CT(\Sigma_{g',1})$ which is convex with respect to the Teichm\"uller metric.
\end{enumerate}
\end{theorem}

The examples constructed in Theorem \ref{non-geometric} are described purely algebraically, but can also be seen as the composition of a monomorphism of the type described by Ivanov-McCarthy, composed with the (noninjective) homomorphism obtained by forgetting some of the marked points. 
\bigskip

We conclude this rather long introduction by a brief plan of the paper. In section \ref{sec:hom cover} we explain the preliminaries on induced homomorphisms from covering spaces, and provide criteria which allow one to construct monomorphisms from covers (see Proposition \ref{liftingcondition}).   Section \ref{sec:main} begins by converting these criteria into statements about surjective homomorphisms of surface groups to finite groups (Proposition \ref{criteria}), and ends by finding homomorphisms satisfying the required properties to prove Theorem \ref{main}.
In Section \ref{sec:alpha} we construct the required monomorphism $\alpha$, and finally in Section \ref{sec:nongeometric} we prove Theorem \ref{non-geometric}.
\bigskip

\noindent{\bf{Acknowledgments.}} The first two authors wish to thank the University of Michigan, where this work began. The first author would like to express his gratitude to the University of Illinois at Urbana-Champaign. The third author would like to thank the Universidad Aut\'onoma de Madrid for its hospitality while this paper was being wirtten. All three of us thank Richard Kent and Alan Reid for, among other things, helpful conversations, and Nathan Dunfield, for pointing out a lemma of Hall (Lemma \ref{hallslem} here) which provides an alternative ending to the construction for Theorem \ref{main}.

\section{Homomorphisms from covers} \label{sec:hom cover}

In this section we describe the basic means of building homomorphisms of mapping class groups from covers.  Maclachlan-Harvey \cite{maclachlanharvey} and Birman-Hilden \cite{Birman-Hilden-Annals} explain this procedure, but because their situation was slightly more complicated---they were allowing branched covers over marked points and then erasing marked points in the covers---the work there was restricted to regular covers.  As our proof of Theorem \ref{main} exploits the irregularity of the covers involved, we cannot refer directly to their work.  Nonetheless, the ideas here are likely well-known, and we have included proofs for completeness.

To keep better track of marked points, we will represent the surface $\Sigma_{g,k}$ by $(\Sigma,{\bf z})$, where $\Sigma$ is a closed surface of genus $g$ and ${\bf z} \subset \Sigma$ is a set of $k$ marked points.  If $k=1$ and ${\bf z} = \{z\}$, we will also write $(\Sigma,z)$.

Let $\kappa: \tilde \Sigma \to \Sigma$ be a degree $d$ cover and ${\bf \tilde z} = \kappa^{-1}({\bf z})$ a set of $kd$ marked points of $\tilde \Sigma$ (thus $(\tilde \Sigma,{\bf \tilde z})$ is a surface of genus $g' = d(g-1)+1$ with $kd$ marked points).  
Given $f \in \Homeo(\Sigma,{\bf z})$ we say that $\tilde f \in \Homeo(\tilde \Sigma,{\bf \tilde z})$ is a lift of $f$ if $f \circ \kappa = \kappa \circ \tilde f$.  If $\tilde f$ exists, then we say that $f$ has a lift.
Define
\[ \Homeo_\kappa(\Sigma,{\bf z})=\{f\in\Homeo(\Sigma,{\bf z}) \, | \, \mbox{there is a lift}\ \tilde f\in\Homeo(\tilde \Sigma,{\bf \tilde z})\ \mbox{of}\  f\}, \]
the group of those homeomorphism of $(\Sigma,{\bf z})$ that have a lift to $(\tilde \Sigma,{\bf \tilde z})$ and
\[\Homeo_\kappa^*(\tilde \Sigma,{\bf \tilde z}) = \{\tilde f \in \Homeo(\tilde \Sigma,{\bf \tilde z}) \, | \, \tilde f \mbox{ is a lift of some}\ f\in\Homeo_\kappa(\Sigma,{\bf z})\},\]
the group of lifts of homeomorphisms of $(\Sigma,{\bf z})$ to $(\tilde \Sigma,{\bf \tilde z})$ by $\kappa$.

There is a homomorphism
\[\kappa_*:\Homeo_\kappa^*(\tilde \Sigma,{\bf \tilde z}) \to \Homeo_\kappa(\Sigma,{\bf z}) \]
with $\kappa_*(\tilde f) = f$ if $\tilde f$ is a lift of $f$.  Two lifts $\tilde f_1$ and $\tilde f_2$ of a given homeomorphism $f$ differ by a covering transformation, and so $\kappa_*$ fits into a short exact sequence
\begin{equation} \label{ses1}
\xymatrix{
1 \ar[r] & K \ar[r] & \Homeo_\kappa^*(\tilde\Sigma,{\bf \tilde z}) \ar[r]^{\kappa_*} & \Homeo_\kappa(\Sigma,{\bf z}) \ar[r] & 1.}
\end{equation}
where $K$ is the deck-transformation group of the covering $\kappa$. We now prove:

\begin{proposition} \label{welldefined}
If $\tilde f \in \Homeo_\kappa^*(\tilde \Sigma,{\bf \tilde z}) \cap \Homeo_0(\tilde \Sigma,{\bf \tilde z})$, then $\kappa_*(\tilde f) = f \in \Homeo_0(\Sigma,{\bf z})$. Moreover, the induced homomorphism
$$\kappa_*:\Homeo_\kappa^*(\tilde \Sigma,{\bf \tilde z}) \cap \Homeo_0(\tilde \Sigma,{\bf \tilde z})\to\Homeo_0(\Sigma,{\bf z})$$
is an isomorphism.
\end{proposition}
\begin{proof}
We begin by replacing both surfaces with punctured surfaces, $\Sigma' = \Sigma - {\bf z}$ and $\tilde \Sigma' = \tilde \Sigma- {\bf \tilde z}$; we still denote the induced cover by $\kappa:\tilde\Sigma'\to\Sigma'$. Fix basepoints $\tilde*\in\tilde\Sigma'$ and $*\in\Sigma'$ with $\kappa(\tilde*)=*$ and consider the associated homomorphism
$$\pi_1(\kappa):\pi_1(\tilde\Sigma',\tilde*)\to\pi_1(\Sigma',*).$$
Choosing a path $\tau$ in $\tilde\Sigma'$ joining the points $\tilde*$ and $\tilde f(\tilde*)$ we obtain an isomorphism between $\pi_1(\tilde\Sigma',\tilde f(\tilde*))$ and $\pi_1(\tilde\Sigma',\tilde*)$. Abusing notation, we then have that $\tilde f$ induces a homomorphism 
$$\pi_1(\tilde f):\pi_1(\tilde\Sigma',\tilde*)\to\pi_1(\tilde\Sigma',\tilde*).$$
Similarly, the projection $\kappa(\tau)$ to $\Sigma'$ of the path $\tau$ we obtain a further homomorphism
$$\pi_1(f):\pi_1(\Sigma',*)\to\pi_1(\Sigma',*).$$
Taking into account that $\tilde f$ is a lift of $f$ we obtain that both homomorphisms $\pi_1(f)$ and $\pi_1(\tilde f)$ satisfy
\begin{equation}\label{eq:welldefined1}
\pi_1(f)\circ\pi_1(\kappa)=\pi_1(\kappa)\circ\pi_1(\tilde f).
\end{equation}
The assumption that $\tilde f\in\Homeo_0(\tilde\Sigma,{\bf\tilde z})$ implies that $\pi_1(\tilde f)$ is an inner automorphism, meaning that there is $\gamma\in\pi_1(\tilde\Sigma',\tilde*)$ with
$$\pi_1(\tilde f)(\eta)=\gamma\eta\gamma^{-1}.$$
Identifying $\pi_1(\tilde\Sigma',\tilde*)$ with its image under $\pi_1(\kappa)$ we obtain thus from \eqref{eq:welldefined1} that the automorphism
\begin{equation}\label{eq:welldefined2}
\pi_1(\Sigma',*)\to\pi_1(\Sigma',*),\ \eta\mapsto\gamma^{-1}(\pi_1(f)(\eta))\gamma
\end{equation}
is the identity on the finite index subgroup $\pi_1(\kappa)(\pi_1(\tilde\Sigma',\tilde*))$ of $\pi_1(\Sigma',*)$. The uniqueness of roots in the surface group $\pi_1(\Sigma',*)$ implies that the homomorphism \eqref{eq:welldefined2} is actually the identity; in other words, $\pi_1(f)$ is an inner automorphism. On the other hand, the automorphism $\pi_1(f)$ preserves the kernel of the homomorphism $\pi_1(\Sigma',*)\to\pi_1(\Sigma,*)$. It follows now from the Dehn-Nielsen-Baer Theorem (see, e.g.~\cite{Farb-Margalit}) that the homeomorphism $f$ is isotopic to the identity. This proves the first claim.

So far, we have shown that the homomorphism
\begin{equation}\label{eq:welldefined3}
\kappa_*:\Homeo_\kappa^*(\tilde \Sigma,{\bf \tilde z}) \cap \Homeo_0(\tilde \Sigma,{\bf \tilde z})\to\Homeo_0(\Sigma,{\bf z})
\end{equation}
is well-defined. The surjectivity of this homomorphism follows directly from the isotopy lifting property of covers. The injectivity can be for instance seen as follows. Endow $\Sigma'$ and $\tilde\Sigma'$ with hyperbolic metrics to that $\kappa$ is a Riemannian cover. Any element $\tilde f$ in the kernel of \eqref{eq:welldefined3} is a lift of the identity and hence an isometry of $\tilde\Sigma'$. It is well-known that the identity is the only isometry of a hyperbolic surface (with non-abelian fundamental group) which is homotopic to the identity. This proves that $\tilde f$ is the identity concluding the proof of Proposition \ref{welldefined}.
\end{proof}

Denote by 
\[\Mod_\kappa^*(\tilde \Sigma,{\bf \tilde z}) = \{ [\tilde f] \in \Mod(\tilde \Sigma,{\bf \tilde z}) \, | \, \tilde f \in \Homeo_\kappa^*(\tilde \Sigma,{\bf \tilde z})\}\]
and
\[\Mod_\kappa(\Sigma,{\bf z}) = \{ [f] \in \Mod(\Sigma,{\bf z}) \, | \, f \in \Homeo_\kappa(\Sigma,{\bf z}) \} \]
the associated subgroups of the mapping class groups of $\Sigma$ and $\tilde\Sigma$, respectively. According to the first claim of Proposition \ref{welldefined}, we can define $\kappa_*:\Mod_\kappa^*(\tilde \Sigma,{\bf \tilde z}) \to \Mod_\kappa(\Sigma,{\bf z})$ by the formula $\kappa_*([\tilde f]) = [\kappa_*(\tilde f)]$. The following is a straightforward consequence of the second claim:

\begin{corollary}\label{cor:well-defined}
The sequence \eqref{ses1} descends to a short exact sequence.
\begin{equation} \label{ses2} \xymatrix{
1 \ar[r] & K \ar[r] & \Mod_\kappa^*(\tilde \Sigma,{\bf \tilde z}) \ar[r]^{\kappa_*} & \Mod_\kappa(\Sigma,{\bf z}) \ar[r] & 1.}
\end{equation}
\end{corollary}

\begin{remark}
Without going any further we would like to remark that the proofs of Proposition \ref{welldefined} and Corollary \ref{cor:well-defined} still apply if $\kappa:(\tilde\Sigma,{\bf\tilde z})\to(\Sigma,{\bf z})$ is a branched cover with all branching points in ${\bf\tilde z}$.
\end{remark}

The groups $\Homeo_\kappa(\Sigma,{\bf z})$ and $\Mod_\kappa(\Sigma,{\bf z})$ can be alternatively described as follows.  Let $\pi_1(\Sigma)$ and $\pi_1(\tilde \Sigma)$ be the fundamental groups of the closed surfaces (the marked points play no role here); as in the proof of Proposition \ref{welldefined} we identify $\pi_1(\tilde\Sigma)$ with a finite index subgroup of $\pi_1(\Sigma)$. The whole group $\Homeo(\Sigma,{\bf z})$ acts on the set of conjugacy classes of finite index subgroups of $\pi_1(\Sigma)$. Since a homeomorphism $f$ has a lift if and only if it preserves the conjugacy class of $\pi_1(\tilde\Sigma)$ in $\pi_1(\Sigma)$, we see that $\Homeo_\kappa(\Sigma,{\bf z})$ is the stabilizer of this conjugacy class.  Similarly, this action of $\Homeo(\Sigma,{\bf z})$ descends to an action of $\Mod(\Sigma,{\bf z})$ on the set of conjugacy classes of finite index subgroups and $\Mod_\kappa(\Sigma,{\bf z})$ is the stabilizer of the conjugacy class of $\pi_1(\tilde\Sigma)$.

The following is now straightforward from the discussion above.
\begin{proposition} \label{liftingcondition}
Given a finite covering $\kappa:\tilde \Sigma \to \Sigma$ identify $\pi_1(\tilde\Sigma)$ with a subgroup of $\pi_1(\Sigma)$. There is an injective homomorphism
\[\Mod(\Sigma,{\bf z}) \to \Mod(\tilde \Sigma,{\bf \tilde z}) \]
obtained by lifting mapping classes in $\Mod(\Sigma,{\bf z})$ to $\Mod(\tilde \Sigma,{\bf \tilde z})$, provided the following two conditions are satisfied
\begin{enumerate}
\item[(i)] the conjugacy class of $\pi_1(\tilde\Sigma)$ in $\pi_1(\Sigma)$ is invariant by the action of $\Mod(\Sigma,{\bf z})$,
\item[(ii)] the sequence \eqref{ses2} is split.
\end{enumerate}
\end{proposition}

We note that condition (ii) can be replaced by requiring \eqref{ses1} to split since this implies \eqref{ses2} is split.

\begin{proof}
By the alternate description of $\Mod_\kappa(\Sigma,{\bf z})$, we see that condition (i) is equivalent to saying that $\Mod(\Sigma,{\bf z}) = \Mod_\kappa(\Sigma,{\bf z})$.  Therefore, the sequence \eqref{ses2} becomes
\begin{equation}\label{ses3} \xymatrix{
1 \ar[r] & K \ar[r] & \Mod_\kappa^*(\tilde \Sigma,{\bf \tilde z}) \ar[r]^{\kappa_*} & \Mod(\Sigma,{\bf z}) \ar[r] & 1.}
\end{equation}
By (ii), \eqref{ses3} splits providing the required injective homomorphism.
\end{proof}

To illustrate the use of Proposition \ref{liftingcondition}, and since we will wish to revisit it later, we now recall the construction due to Ivanov-McCarthy \cite{ivanovmccarthy}.

Let $\Sigma$ have genus $g \geq 2$ and ${\bf z} = \{z\}$, a single point, so that $(\Sigma,z)$ represents the surface $\Sigma_{g,1}$.  Let $\kappa:\tilde \Sigma \to \Sigma$ be a degree $d$ characteristic cover.  By this we mean that $\pi_1(\tilde\Sigma)$ is a characteristic subgroup of $\pi_1(\Sigma)$.  In particular, the conjugacy class of $\pi_1(\tilde\Sigma)$ consists of a single subgroup which is then clearly fixed by the action of $\Mod(\Sigma)$ and so condition (i) of Proposition \ref{liftingcondition} is satisfied.

Let ${\bf \tilde z} = \kappa^{-1}(z)$ be the $d$ marked points in $\tilde \Sigma$, and fix one of them $\tilde z \in {\bf \tilde z}$.
As the cover is characteristic, the covering group $K$ acts transitively on ${\bf \tilde z}$.  So given any $f \in \Homeo(\Sigma,z)$, by composing with a covering transformation if necessary, we can find a lift $\tilde f$ which fixes $\tilde z$.  This specifies a homomorphism $\Homeo(\Sigma,z) \to \Homeo_\kappa^*(\tilde \Sigma,{\bf \tilde z})$ which splits the sequence \eqref{ses1}, and hence condition (ii) of Proposition \ref{liftingcondition} is satisfied.  Therefore, $\kappa$ induces an injective homomorphism 
$$\Mod(\Sigma,z) \to \Mod(\tilde \Sigma,{\bf \tilde z}).$$
In the next section we will make use of a similar strategy to prove Theorem \ref{main}; our covers are however going to be as non-characteristic as possible.

\section{Closed surfaces}\label{sec:main}

In what follows, $\Sigma$ is a closed surface of genus $g \geq 2$ and the marked point set ${\bf z}$ is empty.

\begin{proposition} \label{criteria}
Let $\rho:\pi_1(\Sigma)\rightarrow G$ be a surjective homomorphism to a finite group $G$ such that $\ker(\rho)$ is characteristic. Suppose that $H\subset G$ is a subgroup with
$$(a)\ \ N_G(H)=H\ \ \ \ \ (b)\ \ \Aut(G) \cdot H = \Inn(G) \cdot H.$$
Here $N_G(H)$ is the normalizer of $H$ in $G$ and $\Aut(G)$ and $\Inn(G)$ are the groups of automorphisms and interior automorphisms of $G$, respectively.

Let $\kappa:\tilde{\Sigma} \to \Sigma$ be the cover corresponding to $\rho^{-1}(H)$. Then conditions (i) and (ii) of Proposition \ref{liftingcondition} are satisfied.  In particular, $\kappa$ induces an injective homomorphism
$$\Mod(\Sigma) \to \Mod(\tilde \Sigma).$$
\end{proposition}

\begin{proof}
In order to relax notation we set $\Gamma=\pi_1(\Sigma)$ and $\Gamma_0=\pi_1(\tilde\Sigma)=\rho^{-1}(H)$. The covering group $K$ is isomorphic to the quotient $N_\Gamma(\Gamma_0)/\Gamma_0$ of the normalizer of $\Gamma_0$ in $\Gamma$ modulo $\Gamma_0$. Property (a) implies $N_\Gamma(\Gamma_0) = \Gamma_0$ and hence that $K$ is trivial.  Therefore $\kappa_*$ in \eqref{ses2} is an isomorphism so the sequence splits and (ii) of Proposition \ref{liftingcondition} is satisfied.

To verify condition (i) it suffices to show that the subgroup $\Gamma_0$ has the following property:
\begin{equation}\label{eq:orbits}
\Aut(\Gamma) \cdot \Gamma_0 = \Inn (\Gamma) \cdot \Gamma_0.
\end{equation}
To see that $\Gamma_0$ does indeed have this property, fix $\sigma\in\Aut(\Gamma)$ and observe that since $\rho$ is surjective and has characteristic kernel $\ker(\rho)$ there is $\tau\in\Aut(G)$ with
$$\rho\circ\sigma=\tau\circ\rho.$$
By property (b), there is $g\in G$ such that $\tau(H)=gHg^{-1}$. Since $\rho$ is surjective we have in turn $\gamma\in\Gamma$ with $\rho(\gamma)=g$. We have then
\begin{align*}
\sigma(\Gamma_0)
    &=\sigma(\rho^{-1}(H))=(\rho\circ\sigma^{-1})^{-1}(H)=(\tau^{-1}\circ\rho)^{-1}(H)\\
    &=\rho^{-1}(\tau(H))=\rho^{-1}(gHg^{-1})=\gamma\rho^{-1}(H)\gamma^{-1}=\gamma\Gamma_0\gamma^{-1}.
\end{align*}
This proves \eqref{eq:orbits} and thus verifies condition (i) of Proposition \ref{liftingcondition}.
\end{proof}

To prove Theorem \ref{main}, we must find a group $G$ with a proper subgroup $H < G$ and a surjective homomorphism $\rho:\pi_1(\Sigma)\to G$  as in Proposition \ref{criteria}.  To construct such a group, we consider the $k$-fold product $\CS= S_3 \times \ldots \times S_3$ of $S_3$, the symmetric group in 3 letters, with itself. Let $r_j:\CS \rightarrow S_3$ be the projection onto the $j^{th}$ factor.

\begin{proposition} \label{sylow}
Let $G$ be a subgroup of $\CS$ such that $r_j(G) = S_3$ for all $j = 1,..,k$. If $H$ is a 2-Sylow subgroup of $G$ then $H$ is a proper subgroup satisfying properties (a) and (b) from Proposition \ref{criteria}.
\end{proposition}

\begin{proof}
Since $r_j$ is surjective and $|S_3| = 6$ then $6$ divides $|G|$ and thus any $2$-Sylow subgroup of $G$ is proper. Property (b) is immediate from the fact that all $p$-Sylow subgroups of a finite group are conjugate.

Every $p$-Sylow subgroup of $G$ is the intersection of a $p$-Sylow subgroup of $\CS$ with $G$.  This is because every
$p$-subgroup of $\CS$---in particular a $p$-Sylow subgroup of $G$---is contained in some $p$-Sylow subgroup of $\CS$.

A $2$-Sylow subgroup $P < \CS$ has the following form: there exist $X_1, \ldots ,X_k \in S_3$, each of order 2, such
that
\[ P = \{ \, (x_1, \ldots ,x_k) \, | \, x_j = X_j \mbox{ or } x_j = {\bf 1} \, \} = \langle X_1 \rangle \times  \ldots  \times \langle X_k \rangle, \]
where ${\bf 1}$ is the identity in $S_3$. Let $P$ be such a $2$-Sylow subgroup with $H = P \cap G$.

Now suppose $y = (y_1, \ldots ,y_k) \in N_G(H)$.  We want to show that $y \in H$.  Since $y$ is assumed to lie in $G$,
it suffices to show that $y \in P$.

For every  $j = 1, \ldots ,k$ we claim there exists $x= (x_1, \ldots ,x_k) \in H$ so that $x_j = X_j$. To see this, note that $r_j$ is surjective and so there exists an element $x' \in G$ of order 2 such that $r_j(x')$ has order $2$ in $S_3$. We can conjugate $x'$ by an element of $G$ to some  $x =(x_1, \ldots, x_k) \in H$.  Since $x_j= r_j(x)$ is non-trivial then $x_j = X_j$.

Conjugating $x$ by $y$ we have
\[ ^yx=( ^{y_1}x_1, \ldots , ^{y_k}x_k ) \in H. \]
Therefore, we have
\[ ^{y_j}X_j = X_j \mbox{ or } ^{y_j}X_j = {\bf 1}. \]
The second case is absurd since $X_j \neq {\bf 1}$, so $^{y_j}X_j = X_j$.  Thus $y_j$ is in the centralizer of $X_j$ in
$S_3$ which is $\{{\bf 1},X_j\}$.  Since this is true for all $j$, we have $y_j = X_j$ or $y_j = {\bf 1}$ for all $j$. Therefore,
$y \in P$ as required.
\end{proof}

We can now prove Theorem \ref{main}.

\begin{named}{Theorem \ref{main}}
For every $g\ge 2$ there is $g'>g$ and an injective homomorphism
\[\phi:\Mod(\Sigma_{g,0}) \to \Mod(\Sigma_{g',0}). \]
\end{named}
\begin{proof}
The automorphism group $\Aut(\Gamma)$ of $\Gamma=\pi_1(\Sigma_{g,0})$ acts on $\Hom(\Gamma,S_3)$.  Let $\{\rho_1,\ldots,\rho_k\}$ be the
$\Aut(\Gamma)$--orbit of a surjective homomorphism in $\Hom(\Gamma,S_3)$ (such an epimorphism exists since $S_3$ is 2-generator and $\Gamma$ surjects a rank $2$-free group).  Then
\[ \rho = \rho_1 \times \ldots \times \rho_k:\Gamma \to \CS = S_3 \times \ldots \times S_3 \]
is a homomorphism onto a subgroup $G < \CS$ as in the previous section.  By construction, $\rho$ has characteristic
kernel.  Applying Proposition \ref{criteria} completes the proof.
\end{proof}

We will now describe an alternative construction, appealing again to Proposition \ref{criteria}, using a result of P.~Hall \cite{Hall} about finite simple groups which was explained to us by N.~Dunfield (see also Lemma 3.7 of Dunfield--Thurston \cite{dunthur}).  One version of Hall's Theorem is as follows.

\begin{lemma} [Hall] \label{hallslem}
Suppose $\Gamma$ is any group, that $Q_1, \ldots , Q_k$ are finite, non-abelian simple groups, and that $\rho_i:G \to Q_i$ is an epimorphism for $i = 1,\ldots,k$.
If $\rho_i$ and $\rho_j$ do not differ by an isomorphism $Q_i \to Q_j$ for any $i \neq j$, then
\[ \rho= \rho_1 \times \ldots \times \rho_k: \Gamma \to Q_1 \times \ldots \times Q_k \]
is surjective.
\end{lemma}

Now consider, for example, a finite simple group of the form $\PSL_2(\mathbb F_p)$ (this group is simple when $p$ is a prime greater than or equal to $5$).
Let $\rho_0: \Gamma = \pi_1(\Sigma) \to \PSL_2(\mathbb F_p)$ be an epimorphism.  Again, $\Aut(\Gamma)$ acts on $\Hom(\Gamma,\PSL_2(\mathbb F_p))$, and we note that $\Aut(\PSL_2(\mathbb F_p))$ also acts by postcomposition.  Let $\{\rho_1,...,\rho_n\}$ be a maximal collection of elements of an $\Aut(\Gamma)$--orbit, no two of which are in the same $\Aut(\PSL_2(\mathbb F_p))$--orbit.  This is easily accomplished by starting with any $\Aut(\Gamma)$--orbit and discarding some number of elements.

The point is that now $\rho_i$ and $\rho_j$ do not differ by an automorphism of $\PSL_2(\mathbb F_p)$ (unless $i = j$), and by maximality, given $\phi \in \Aut(\Gamma)$ and $i \in \{ 1,...,n\}$, there exists $\tau \in \Aut(\PSL_2(\mathbb F_p))$ and $j \in \{1,...,n\}$ so that $\rho_i \circ \phi = \tau \circ \rho_j$.  Therefore, if we set
\[ \rho = \rho_1 \times \ldots \times \rho_k: \Gamma \to \PSL_2(\mathbb F_p) \times \ldots \times \PSL_2(\mathbb F_p) = G,\]
then $\rho$ has characteristic kernel and by Lemma \ref{hallslem}, $\rho$ is surjective.

Now, let $H_0 < \PSL_2(\mathbb F_p)$ to be the subgroup of upper triangular matrices and set $H = H_0 \times \ldots \times H_0 < G$.  Since $N_{\PSL_2(\mathbb F_p)}(H_0) = H_0$, it follows that $N_G(H) = H$ so condition (a) from Proposition \ref{criteria} is satisfied.  Furthermore, as with any product of nonabelian finite simple groups, $\Aut(G)$ acts by permuting the factors and composing with automorphisms in each factor.  The subgroup $\Aut(\PSL_2(\mathbb F_p))$ contains $\Inn(\PSL_2(\mathbb F_p))$ as an index two subgroup, and it is not hard to see that $\Aut(\PSL_2(\mathbb F_p)) \cdot H_0 = \Inn(\PSL_2(\mathbb F_p)) \cdot H_0$.  It follows that $\Aut(G) \cdot H = \Inn(G) \cdot H$, and so condition (b) from Proposition \ref{criteria} is also satisifed.

\bigskip

At this point we would like to observe that the proof of Theorem \ref{main} also applies in other situations. For instance, we obtain injective homeomorphisms between the mapping class groups of two non-homeomorphic handlebodies. Similarly, we obtain injective homomorphisms
\begin{equation}\label{outfn}
\Out(\BF_n)\to\Out(\BF_m)
\end{equation}
with $n<m$ where $\Out(\BF_n)$ is the group of outer automorphisms of the free group $\BF_n$. The existence of such homomorphisms \eqref{outfn} was previously obtained by Bogopol'skii-Puga \cite{Bogopolskii-Puga} using a different argument. In the examples of \cite{Bogopolskii-Puga}, $m$ grows exponentially with $n$. However, the authors have only been able to obtain doubly exponential bounds on the degree of the covers constructed in the proof of Theorem \ref{main}. Thus we ask:

\begin{question}
What is the minimal degree covering $\tilde \Sigma \to \Sigma$ of a closed surface $\Sigma$ by $\tilde \Sigma$ for which one can find an injection $\Mod(\Sigma)\to\Mod(\tilde \Sigma)$?
\end{question}

\section{Another injective homomorphism}\label{sec:alpha}

We now set out to construct the homomorphism $\alpha$ of Theorem \ref{non-geometric}.

We assume that $\Sigma$ and $\tilde \Sigma$ are closed surfaces of genus $g$ and $g'$, respectively, and that $\kappa:\tilde \Sigma \to \Sigma$ is a characteristic covering.  Let $z \in \Sigma$ and $\tilde z \in \tilde \Sigma$ be single marked points with $\kappa(\tilde z) = z$ and let $\Gamma = \pi_1(\Sigma,z)$ and $\tilde \Gamma = \pi_1(\tilde \Sigma,\tilde z)$.  As before, we write $\pi_1(\kappa): \tilde \Gamma \to \Gamma$ for the induced homomorphism on fundamental groups. We write $(\Sigma,z)$ and $(\tilde\Sigma,\tilde z)$ for the surfaces with marked points.  So we wish to construct a homomorphism
\[\alpha:\Mod(\Sigma,z) \to \Mod(\tilde \Sigma,\tilde z). \]

The Birman exact sequence \cite{Birman}, is given by
\[ \xymatrix{ 1 \ar[r] & \Gamma \ar[r] & \Mod(\Sigma,z) \ar[r] & \Mod(\Sigma) \ar[r] & 1.}\]
An element of $\Mod(\Sigma,z)$ induces an isomorphism of $\Gamma$ and an element of $\Mod(\Sigma)$ induces an outer automorphism.  Since an element of $\Gamma$ clearly gives us an inner automorphism, one can check that this determines a homomorphism of short exact sequences
\begin{equation} \label{sespair} \xymatrix{ 1 \ar[r] & \Gamma \ar[r] \ar[d] & \Mod(\Sigma,z) \ar[r] \ar[d] & \Mod(\Sigma) \ar[r] \ar[d] & 1 \\
1 \ar[r]  & \Inn(\Gamma) \ar[r] & \Aut(\Gamma) \ar[r] & \Out(\Gamma) \ar[r] & 1. }
\end{equation}
Since $\Gamma$ has trivial center, the first vertical arrow is an isomorphism.  Moreover, by the Dehn-Nielsen-Baer Theorem, the second and third vertical arrows are isomorphisms. A similar discussion holds for $\tilde \Sigma$ and $\tilde \Gamma$.  We will use the central vertical homomorphisms to identify $\Mod(\Sigma,z) = \Aut(\Gamma)$ and $\Mod(\tilde \Sigma,\tilde z) = \Aut(\tilde \Gamma)$.  Similarly, we use the initial vertical arrows to identify $\Inn(\Gamma) = \Gamma$ and $\Inn(\tilde \Gamma) = \tilde \Gamma$.

We are now required to construct
\[ \alpha:\Aut(\Gamma) \to \Aut(\tilde \Gamma). \]
Given $\phi \in \Aut(\Gamma)$, we define $\alpha(\phi) \in \Aut(\tilde \Gamma)$ by the formula
\[\alpha(\phi)(\tilde \gamma) =\pi_1(\kappa)^{-1}(\phi(\pi_1(\kappa)\tilde \gamma)) \]
for every $\tilde \gamma \in \tilde \Gamma$. This makes sense because $\pi_1(\kappa)$ is an isomorphism onto its image, $\pi_1(\kappa)(\tilde \Gamma)$, which is a characteristic subgroup of $\Gamma$.  Thus, we are simply restricting $\phi$ to $\pi_1(\kappa)(\tilde \Gamma)$, and conjugating it back to $\tilde \Gamma$ by $\pi_1(\kappa)$.

\begin{lemma}\label{lem2}
The homomorphism $\alpha:\Aut(\Gamma) \to \Aut(\tilde \Gamma)$ is injective.
\end{lemma}
\begin{proof}
Given $\phi \in \Aut(\Gamma)$, we assume that $\alpha(\phi)$ is the identity.  So, $\alpha(\phi)(\tilde \gamma) = \tilde \gamma$ for all $\tilde \gamma \in \tilde \Gamma$, and so $\phi(\gamma) = \gamma$ for all $\gamma \in\pi_1(\kappa)(\tilde \Gamma)$.  Since $\pi_1(\kappa)(\tilde \Gamma)$ has finite index in $\Gamma$ and $\Gamma$ has unique roots, it follows that $\phi(\gamma) = \gamma$ for all $\gamma \in \Gamma$.  That is, $\phi$ is the identity, and $\ker(\alpha)$ is trivial as required.
\end{proof}

One can alternatively describe $\alpha$ as follows.  The Ivanov-McCarthy monomorphism $\Mod(\Sigma,z) \to \Mod(\tilde \Sigma,{\bf \tilde z})$ as described in Section \ref{sec:hom cover} has image inside $\Mod(\tilde \Sigma,{\bf \tilde z},\tilde z)$, the group of those homeomorphisms of $(\tilde \Sigma,{\bf \tilde z})$ which fix the single marked point $\tilde z$.  There is a homomorphism (fitting into Birman's more general exact sequence \cite{Birman}) defined by forgetting ${\bf \tilde z} - \tilde z$, all the marked points \textit{except} $\tilde z$
\[\Mod(\tilde \Sigma,{\bf \tilde z},\tilde z) \to \Mod(\tilde \Sigma,\tilde z).\]
We leave it as an exercise for the interested reader to check that $\alpha$ is given as the composition
\[\Mod(\Sigma,z) \to \Mod(\tilde \Sigma,{\bf \tilde z},\tilde z) \to \Mod(\tilde \Sigma,\tilde z). \]

\section{Properties of $\alpha$} \label{sec:nongeometric}

In this section we prove Theorem \ref{non-geometric} but before doing so we motivate it briefly. The examples of injective homomorphisms between different mapping class groups constructed in \cite{Birman-Hilden}, \cite{ivanovmccarthy} or Theorem \ref{main} are induced by a cover in the following sense.

\begin{definit}
Let $(\Sigma,{\bf z})$ and $(\tilde\Sigma,{\bf\tilde z})$ be surfaces with a possibly empty finite set of points marked. A homomorphism $\phi:\Mod(\Sigma,{\bf z})\to\Mod(\tilde\Sigma,{\bf\tilde z})$ is {\em induced by a cover} if there is a possibly branched cover $\kappa:\tilde\Sigma\to\tilde\Sigma$, with all branching points of $\kappa$ contained in ${\bf z}$, with ${\bf\tilde z}\subset\kappa^{-1}{\bf z}$, such that if $z_0\in\kappa^{-1}({\bf z})\setminus {\bf \tilde z}$ then $\deg_\kappa(z_0)>1$ and such that there is a continuous homomorphism 
$$\kappa_*:\Homeo(\Sigma,{\bf z})\to\Homeo(\tilde\Sigma,{\bf\tilde z}),$$
for which $\kappa_*(f)$ is a lift of $f$ for all $f \in \Homeo(\Sigma,{\bf z})$, and which induces $\phi$.
\end{definit}

If $\phi:\Mod(\Sigma,{\bf z})\to\Mod(\tilde\Sigma,{\bf\tilde z})$ is induced by a cover $\kappa:(\tilde\Sigma,{\bf\tilde z})\to(\Sigma,{\bf z})$ such that ${\bf \tilde z}=\kappa^{-1}({\bf z})$ then it follows from Proposition \ref{welldefined} and the remark following Corollary \ref{cor:well-defined}, that there is a well-defined $\phi$-equivariant isometric embedding
$$\Phi:\CT(\Sigma,{\bf z})\to\CT(\tilde\Sigma,{\bf\tilde z}).$$
The result remains true if ${\bf\tilde z}$ is properly contained in $\kappa^{-1}({\bf z})$ since the assumption on the local degree of the cover implies that the composition of the isometric embedding $\CT(\Sigma,{\bf z})\to\CT(\tilde\Sigma,\kappa^{-1}({\bf z}))$ with the projection $\CT(\tilde\Sigma,\kappa^{-1}({\bf z}))\to\CT(\tilde\Sigma,{\bf\tilde z})$ is still an isometric embedding. In other words, we have:


\begin{lemma}\label{geometric}
Assume that a homomorphism $\phi:\Mod(\Sigma,{\bf z})\to\Mod(\tilde\Sigma,{\bf\tilde z})$ is induced by a cover and endow the Teichm\"uller spaces $\CT(\Sigma,{\bf z})$ and $\CT(\tilde\Sigma,{\bf\tilde z})$ with the Teichm\"uller metric. Then there is an $\phi$-equivariant isometric embedding
$$\Phi:\CT(\Sigma,{\bf z})\to\CT(\tilde\Sigma,{\bf\tilde z}).$$
In particular $\phi$ is type-preserving in the sense that it maps pseudo-Anosov and reducible elements of $\Mod(\Sigma,{\bf z})$ to again pseudo-Anosov and reducible elements of $\Mod(\tilde\Sigma,{\bf\tilde z})$. Also, the limit set $\Lambda_{\phi(\Mod(\Sigma,{\bf z}))}$ of the image of $\phi$ is a proper subset of $\PML(\tilde\Sigma,{\bf\tilde z})$ homeomorphic to $\PML(\Sigma,{\bf z})$. \qed
\end{lemma}

Here we have identified $\PML(\tilde\Sigma,{\bf\tilde z})$ with the Thurston boundary of the Teichm\"uller space $\CT(\tilde\Sigma,{\bf\tilde z})$. McCarthy-Papadopoulos \cite{McCarthy-Papadopoulos} defined the limit set $\Lambda_H$ of a subgroup $H\subset\Mod(\tilde \Sigma,\tilde z)$ to be the closure in $\PML(\tilde\Sigma,\tilde z)$ of the set of pseudo-Anosov fixed points (when there are no pseudo-Anosov elements, the definitions must be modified but this is not relevant to us since, as we just remarked, $\phi(\Mod(\Sigma,{\bf z}))$ contains such elements). 

In a nutshell, the claim of Theorem \ref{non-geometric} is that none of these properties holds for the homomorphism $\alpha$ constructed in the last section.

\begin{named}{Theorem \ref{non-geometric}}
For every $g$ there is $g'$ and an injective homomorphism $\alpha:\Mod(\Sigma_{g,1})\to\Mod(\Sigma_{g',1})$ with the following properties:
\begin{enumerate}
\item There are pseudo-Anosov mapping classes in $\Mod(\Sigma_{g,1})$ whose image under $\alpha$ is a multi-twist. In particular, $\alpha$ is not type-preserving.
\item The limit set $\Lambda_{\alpha(\Mod(\Sigma_{g,1}))}$ of the image of $\alpha$ is the whole space $\PML(\Sigma_{g',1})$ of projective measured laminations on $\Sigma_{g',1}$. 
\item In particular, there is no $\alpha(\Mod(\Sigma_{g,1}))$-invariant subset of $\CT(\Sigma_{g',1})$ which is convex with respect to the Teichm\"uller metric.
\end{enumerate}
\end{named}

With the same notation as in the preceding section, we use the Birman exact sequence \eqref{sespair} to view $\Gamma=\pi_1(\Sigma,z)$ and $\tilde\Gamma=\pi_1(\tilde\Sigma,\tilde z)$ as subgroups of $\Mod(\Sigma,z)$ and $\Mod(\tilde \Sigma,\tilde z)$ respectively. The proof of Theorem \ref{non-geometric} is based on the following observation.

\begin{lemma} \label{pre-denselimitset}
$\alpha(\Gamma)$ contains $\tilde \Gamma$ with finite index.
\end{lemma}
\begin{proof}
From the definition we have $\alpha(\pi_1(\kappa)(\tilde \Gamma)) = \tilde \Gamma$, and $\pi_1(\kappa)(\tilde \Gamma) < \Gamma$ has finite index.
\end{proof}

We start with the proof of the first claim of Theorem \ref{non-geometric}. To begin with, recall that in \cite{kra}, Kra proved that $\gamma \in \Gamma$ is pseudo-Anosov if and only if the free homotopy class determined by $\gamma$ is filling in $\Sigma$ --that is, every element of the free homotopy class of $\gamma$ nontrivially intersects every homotopically essential closed curve. On the other hand, if $\gamma$ is freely homotopic to a simple curve, then the mapping class determined by $\gamma$ is a multi-twist. In particular, in order to find elements $[f]\in\Mod(\Sigma_{g,1})$ of the required type, it suffices to exhibit an element $\gamma\in\Gamma$ whose free homotopy class is filling, but which lifts to a simple loop on $\tilde \Sigma$.  Such loops $\gamma$ can be constructed as follows. Let $\lambda$ be an ending lamination (meaning a measurable geodesic lamination with all complementary regions being ideal polygons) on $\tilde \Sigma$ which is not the preimage of a lamination on $\Sigma$.  Let $\beta_n$ be a sequence of simple closed geodesics in $\Sigma$ which are converging to $\lambda$ in the Hausdorff topology.  Since $\lambda$ is not a lift, we can find an element $\delta \in K$, the covering group of $\kappa:\tilde \Sigma \to \Sigma$, for which $\delta(\lambda) \neq \lambda$.  Since $\lambda$ is an ending lamination, $\lambda$ and $\delta(\lambda)$ fill up $\tilde \Sigma$.  In particular, for $n$ sufficiently large, $\beta_n$ and $\delta(\beta_n)$ also fill up $\tilde \Sigma$.   It follows that for sufficiently large $n$, $\kappa(\beta_n)$ is filling on $\Sigma$ and has a simple lift, $\beta_n$, to $\tilde \Sigma$.  Appropriately replacing the curves by based loops gives examples of the required type for all sufficiently large $n$. This concludes the proof of the first claim of Theorem \ref{non-geometric}.

\begin{remark}
Obviously, if $\gamma\in\tilde\Gamma$ fills $\tilde\Sigma$, then its projection also fills $\Sigma$. In particular, there are pseudo-Anosov mapping classes in $\Mod(\Sigma,z)$ whose image under $\alpha$ is also pseudo-Anosov. 
\end{remark}

We now prove the second claim of Theorem \ref{non-geometric}. 
If $H\subset\Mod(\tilde\Sigma,{\bf\tilde z})$ is a subgroup, let $\Lambda_H$ be its limit set in the sense of McCarthy-Papadopoulos \cite{McCarthy-Papadopoulos}; we will only consider subgroups containing pseudo-Anosov elements. By Lemma \ref{pre-denselimitset}, the image $\alpha(\Gamma)$ of $\Gamma$ contains $\tilde \Gamma$ with finite index and hence
$$\Lambda_{\tilde\Gamma}\subset\Lambda_{\alpha(\Mod(\Sigma,z))}.$$
On the other hand, $\tilde \Gamma$ is an infinite normal subgroup containing two non-commuting pseudo-Anosov elements, and hence it follows from \cite[Prop. 5.5]{McCarthy-Papadopoulos} that
\[ \Lambda_{\tilde \Gamma} =\Lambda_{\Mod(\tilde \Sigma,\tilde z)}=\PML(\tilde\Sigma,\tilde z).\]
The second claim of Theorem \ref{non-geometric} follows directly from the two last equations.

Finally, let us establish the third claim of Theorem \ref{non-geometric}.  It follows from the definition of the limit set and from the dynamics of pseudo-Anosov elements that if $H\subset\Mod(\tilde\Sigma,\tilde z)$ is a subgroup (containing pseudo-Anosovs) then for any $[\lambda] \in \Lambda_H$ and $X \in \CT(\tilde \Sigma,\tilde z)$, there is a sequence $\{g_i\} \subset H$ for which $\lim_ig_iX=[\lambda]$. We make use of the following result:

\begin{lemma}\label{lem3}
Let $S$ be a surface of finite analytical type and assume that a subgroup $H\subset\Mod(S)$ contains pseudo-Anosov elements and stabilizes a proper, closed, convex subset of $\CT(S)$ with respect to the Teichm\"uller metric. Then $\Lambda_H$ is a proper subset of $\PML(S)$.
\end{lemma}

Recall that a subset $Z$ of a metric space $X$ is convex if any geodesic segment in $X$ with endpoints in $Z$ is itself contained in $Z$.

\begin{proof}[Proof of Lemma \ref{lem3}]
Assume that $\Lambda_H=\PML(S)$ and that $\CZ\subset\CT(S)$ is an $H$-invariant convex closed subset; fix a point $X\in\CZ$. Given $[\lambda]\in \PML(S)$ uniquely ergodic, there is \cite{Hubbard-Masur} a unique geodesic ray $\eta^\lambda(t)$ in $\CT(S)$ with $\eta^\lambda(0)=X$ and $\lim_{t\to\infty}\eta^\lambda(t)=[\lambda]$. Since $[\lambda]\in\Lambda_H$ there is a sequence $\{g_i\}\subset H$ with $\lim_ig_iX=[\lambda]$.

Up to subsequence, the geodesics $[X,g_iX]$ converge to a geodesic ray from $X$.  By the unique ergodicity of $[\lambda]$, the limiting ray must be $\eta^\lambda$.
Since the set $\CZ$ is convex and $X$ and $g_iX$ both belong to $\CZ$, we have that $[X,g_iX]\subset\CZ$. By continuity, this implies that the ray $\eta^\lambda\subset\CZ$ as well. We derive that $\CZ$ contains all rays starting in $X$ and tending to an uniquely ergodic point in $\PML(S)$. By \cite{Masur-annals}, this set of rays is dense in $\CT(S)$; hence $\CZ=\CT(S)$.

We have proved that if $\Lambda_H=\PML(\tilde\Sigma,\tilde z)$ then $\CT(S)$ does not contain any proper, closed, non-empty convex $H$-invariant sets.
\end{proof}

It follows directly from the second claim of Theorem \ref{non-geometric} and Lemma \ref{lem3} that there is no $\alpha(\Mod(\Sigma,z))$-invariant subset of $\CT(\tilde\Sigma,\tilde z)$ which is convex with respect to the Teichm\"uller metric. This concludes the proof of Theorem \ref{non-geometric}.\qed
\medskip

Before moving on we would like add some observations. To begin with, observe that the same argument used in \cite[Theorem 2.4]{Kent-Leininger} shows that there is no $\alpha$-equivariant continuous map $\PML(\Sigma,z)\to \PML(\tilde\Sigma,\tilde z)$. Observe also that the arguments given above can be used to prove that there is also no $\alpha(\Mod(\Sigma,z))$-invariant proper subset of $\CT(\tilde\Sigma,\tilde z)$ which is convex with respect to the Weil-Petersson metric. The needed facts about the Weil-Petersson metric can be found in \cite{BMM}; more precisely, see Theorem 1.4 and Theorem 1.8 therein.

From some point of view, it may be more natural to consider the Weil-Petersson than the Teichm\"uller metric. For a start, observe that if a homomorphism is induced by a non-branched cover 
then the map provided by Lemma \ref{geometric} is, up to scaling, also an isometric embedding with respect to the Weil-Petersson metric; its image is totally geodesic. Also, if $\Sigma\subset\Sigma'$ is a compact subsurface then the mapping class group of $\Sigma$ is (usually) isomorphic to a subgroup of the mapping class group of $\Sigma'$, see \cite{parisrolfsen}, \cite{bellmargalit}. We have then an equivariant map from the Teichm\"uller space of $\Sigma$ to the {\em Weil-Petersson completion} of the Teichm\"uller space of $\Sigma'$ whose image is totally geodesic and hence convex.
\medskip

\noindent {\bf Something positive about $\alpha$:} While there are no $\alpha$-equivariant isometric embeddings of $\CT(\Sigma,z)$ to $\CT(\tilde \Sigma,\tilde z)$ it is not difficult to construct such a map which is holomorphic. In fact, the covering $\tilde\kappa$ induces an equivariant holomorphic embedding $\CT(\Sigma,z)\to\CT(\tilde \Sigma,{\bf \tilde z})$. Filling all the punctures but $\tilde z$ yields then a holomorphic fibration $\CT(\tilde \Sigma,{\bf \tilde z})\to\CT(\tilde \Sigma,\tilde z)$, see \cite{bersfiber}. The composition
$$F:\CT(\Sigma,z)\to\CT(\tilde\Sigma,\tilde z)$$
is the desired $\alpha$-equivariant holomorphic map. Observe that being holomorphic, $F$ is $1$-Lipschitz with respect to the Teichm\"uller metric; this follows from the fact that the Teichm\"uller metric and the Kobayashi metric coincide on Teichm\"uller space. 

A brief computation using the Ahlfors Lemma proves that the projection map $\CT(\tilde\Sigma,{\bf\tilde z})\to\CT(\tilde \Sigma,\tilde z)$ is 1-Lipschitz with respect to the Weil-Petersson metrics. On the other hand, the lifting map $\CT(\Sigma,z)\to\CT(\tilde\Sigma,{\bf\tilde z})$ is, again with respect to the Weil-Petersson metrics, a homotethy onto its image. In particular, the map $F$ is Lipschitz with respect to the Weil-Petersson metrics on the source and target Teichm\"uller spaces. This implies for example that the induced map 
$$f:\CM(\Sigma,z)\to\CM(\tilde\Sigma,\tilde z)$$
between the corresponding moduli spaces has finite energy. On the other hand, being holomorphic, the map $f$ is harmonic \cite[8.15]{Eells} and its image is minimal by the Wirtinger inequality.
\medskip

\noindent {\bf Musings on a questions of Farb and Margalit:} As mentioned in the introduction, it is an interesting problem to study how injective homomorphisms between mapping class groups arise. In this spirit, Farb-Margalit (see Question 2 in \cite{bellmargalit}) have asked:

\begin{question}[Farb-Margalit]
Is every injective homomorphism of mapping class groups geometric?
\end{question}

\noindent A problem with this question is that the word {\em geometric} does not have a precise meaning. In some sense, Theorem \ref{non-geometric} asserts that it is difficult to find an interpretation of the word geometric for which the answer can be positive (at least as long as punctures are not explicitly avoided). This is a sad fact since one would hope that the answer to this question is ``yes''. A possible interpretation of the word geometric could be that every such homomorphism induces an algebraic morphism between the corresponding moduli spaces. Another may be simply that it is induced by some manipulation of surfaces. In fact, it may turn out that the second definition is equivalent to the first one.

\bigskip

\noindent Department of Mathematics, National University of Ireland, Galway, Ireland \newline \noindent
\texttt{Javier.Aramayona@nuigalway.ie}

\bigskip

\noindent Department of Mathematics, University of Illinois, Urbana--Champaign, IL 61801 \newline \noindent
\texttt{clein@math.uiuc.edu}

\bigskip

\noindent Department of Mathematics, University of Michigan, Ann Arbor, MI 48109 \newline \noindent
\texttt{jsouto@umich.edu}

\end{document}